\documentstyle[twoside]{article}

\newtheorem{corollary}{Corollary}
\newtheorem{theorem}{Theorem}
\newtheorem{lemma}{Lemma}
\newtheorem{remark}{Remark}

\begin{document}

\begin{center}
{\bf { ON FINITENESS OF CERTAIN VASSILIEV INVARIANTS
}}

\vspace{5mm}

{LOUIS H. KAUFFMAN \footnote{Supported in part by NSF
Grant No. DMS-9504471} }

{ {\it Department of Mathematics}}

 { { \it University of Illinois at Chicago}}

{ {\it 851 S. Morgan St.}}

{ { \it Chicago, Illinois 60607-7045}}

\&

{ MASAHICO SAITO}

{ {\it Department of Mathematics}}

{ {\it Northwestern University }}

{{ \it Evanston, Illinois 60208 }}

\&

{ STEPHEN F. SAWIN \footnote{Supported in part by NSF
Grant \#23068} }

{ {\it Department of Mathematics}}

{ {\it M. I. T.  }}

{{ \it Cambridge, Massachusetts 02139 }}

\end{center}

\rule{0in}{1in}

\begin{abstract}
The best known examples of Vassiliev invariants are the
coefficients of a Jones-type  polynomial expanded after exponential
substitution.
We show that for a given knot, the first $N$ Vassiliev invariants
in this family determine the rest for some integer $N$.
\end{abstract}

\newpage

In this paper we prove a finiteness property for  certain families
of Vassiliev knot invariants.
Let $V_K(t)$ denote the Jones polynomial \cite{Jones}
of a classical knot or a link $K$
(i.e. smoothly embedded circle(s) in the $3$-space).
It is known \cite{BL,B} that the coefficients of powers of $x$
in $V_K(e^x)$ are Vassiliev invariants of finite type.
More specifically, if we write
$V_K(e^x)= \sum_{n=1}^{\infty}v_n(K)x^n$,
then $v_n(K)$ is a Vassiliev invariant of type $n$.
Let us call $\{ v_n(K) \}$ the {\it Jones-Vassiliev } invariants.

More generally, if $P_K(\alpha, z)$ denotes the HOMFLY
polynomial \cite{F,Jones}  and $Y_K(\alpha,z)$ the Kauffman polynomial
\cite{K2},
we obtain a Laurent polynomial (possibly times $t^{1/2}$)
by
substituting
$t^{1/2}-t^{-1/2}$ for $z$ and $t^a$
for $\alpha$, where
$a$ is a nonnegative
half-integer.    Again substituting $e^x$ for
$t$, the $n$th coefficient of the power series in $x$ is a degree $n$
Vassiliev invariant \cite{BL,B}.
Let us call these Vassiliev invariants respectively $p_{a,n}(K)$ and
$f_{a,n}(K)$.   We prove the following.

\begin{theorem}
Given a projection of a  link $K$,  $v_n(K)$ is a linear combination of
the quantities $v_0(K), \ldots ,v_N(K)$, where $N$ is
the number of crossings in the
projection.
The coefficients of the linear combination depend on
the following data: $n$, $N$, the number of
separable pieces in $K$, the writhe of the projection, and the
number of components of the link obtained by parallel smoothings
at positive crossings and transverse smoothings at negative crossings
that are defined in \cite{Mura}.
\end{theorem}

\begin{theorem}
Given a projection of a link $K$ and a half-integer $a$,
$p_{a,n}(K)$
is a linear combination of $p_{a,0}(K), \ldots p_{a,N}(K)$,
where $N$ is
$(2a-1)(s-1) + c$,
$s$ is the number of Seifert circles in the diagram, and
$c$ is the crossing number
of the link.
The coefficients of the linear combination
depend only on $n$, $N$, and the writhe of the diagram.
\end{theorem}

\begin{theorem}
Given a projection of a link $K$ and a half-integer $a$,
$f_{a,n}(K)$
is a linear combination of $f_{a,0}(K), \ldots f_{a,N}(K)$
where $N$ is the greatest integer less than $3ac$, for $c$ the
crossing number, and  the linear combination depend only on
$n$, $N$,  the number of positive and the number of negative
crossings in the diagram, and the number of connected components of
the diagram.
\end{theorem}

These theorems follow from bounds on the degrees of the link
polynomials together with the following straightforward lemma.

\begin{lemma} If $F(t)$ is of the form $t^M$ times a
degree $N$ polynomial $P(t)$ in $t$,
so that $F(t)=t^M P(t)$,  for some $M \in {\bf R}$, %{\mathbb R}$,
and the
largest and smallest powers of $t$ in $F$ are $L$ and $M$
respectively (so that $N=L-M$), then there exist
linear functions $f_{L,N,j}(x_0, \dots, x_{N})$, for each
positive integer $j$, such that
$$d^j/dx^j F(e^x)|_{x=0} = f_{L,N,j}(F(1),d/dx(F(e^x))|_{x=0}, \ldots,
d^{N}/dx^{N}(F(e^x))|_{x=0}).$$
\end{lemma}
%%%%%%%%%%%%%%%%%%%%%%%%Alternate proof%%%%%%%%%%%%%

{\it Proof.}
Consider the linear ODE
$$ D^{N+1} e^{-Mx} f(x)=0,$$
where $D= d/d(e^x)= e^{-x}d/dx$.  Notice that changing variables via
$t=e^x$ this
is just $(d^{N+1}/dt^{N+1}) t^{-M}f(\ln t)=0$, so in particular $f(x)=F(e^x)$
is
a solution.  Since this is a linear $N+1$st order differential
equation in $x$,
 any solution is a linear combination
of the $N+1$ fundamental solutions with the coefficients being
certain linear combinations of the
initial data.
By the initial data we mean
the first $N$ derivatives of the solution
evaluated at $x=0$, which in this case are
$F(1),d/dx(F(e^x))|_0, \ldots,
d^{L-M}/dx^{L-M}(F(e^x))|_0$.
By differentiating this linear combination $j$ times, we see that
$(d^j/dx^j)F(e^x)|_0$ is a linear combination of the initial data too.
$\Box$

\vspace{24pt}

{\it Proof (of Theorem 1).}
By \cite{K3,Mura,Mowen} the span of the Jones polynomial (the highest
degree minus
the lowest degree) is bounded by $c+ g -1$, where $c$ is the
crossing number and $g$ is the number of
disconnected components of the projection.  Further, the upper and
lower bounds depend on the quantities given in the statement of the
theorem.
 But
notice that $V_K(t)$
is divisible by $(t^{1/2}+t^{-1/2})^{g-1}$, since this is true of a
link with $g$ or more unlinked unknots, and $V_K(t)$ is a linear
combination of such by the skein relations.  The result of dividing
$V_K(t)$ by this factor is a Laurent polynomial (possibly times
$t^{1/2}$) with span bounded by exactly $c$.  By the lemma,  the
$j$th derivative of this Laurent polynomial evaluated at $t=1$ can be
written explicitly as a linear combination of the various $k$th
derivatives for $0 \leq k \leq c$.  Since the $j$th derivative of
$V_K(t)$ is a linear combination of the $k$th derivative of this
polynomial for $k \leq j$, and vice versa,  this gives the result.
$\Box$

\begin{corollary}  If a knot $K$ has a projection with $N$ crossings,
  and  $v_k(K)=0$ for $1 \leq k \leq N$, then the
  Jones polynomial
  of $K$ is $1$. $\Box$
\end{corollary}

\vspace{24pt}

{\it Proof (of Theorem 2).}
It was proved in  \cite{Mor} that
$$  d_{\max}(z) \leq c - s + 1$$
and
$$ w - s + 1 \leq d_{\min}(\alpha)
\leq d_{\max}(\alpha) \leq w + s -1 $$
where $d_{\min}$ and $d_{\max}$ denote the lowest and highest degrees of
the respective variable in $P_K(\alpha,z)$, and $c$, s and $w$ are
respectively the crossing number,
number of Seifert circles and writhe
of the projection of a knot $K$.

Notice here that the negative powers of $z$ arise from
the {\it loop value} (the contribution of a disjoint trivial
circle in the skein computations, see \cite{K2})
$ \delta =  (t^a-t^{-a})/(t^{1/2}-t^{-1/2}) $
which is a power of $t$ times a polynomial in $t$.
Thus terms with negative
powers of $z$ have lower $d_{\max}(t)$ and higher $d_{\min}(t)$
 than
the degree of $\alpha$ would indicate.
So the largest power of $t$ which can occur would come from the
product of the largest power of $\alpha$  with the largest power of
$z$, and the smallest power of $t$ comes from the smallest power of
$\alpha$ times the largest power of $z$.
Since
the largest power of $z$ occurs multiplied by the lowest and highest power
of $\alpha$, we have after substitution that
$$  a \cdot d_{\min}(\alpha) - d_{\max}(z)/2 \leq d_{\min} (t) \leq
d_{\max}(t) \leq a \cdot d_{\max}(\alpha) + d_{\max}(z)/2  $$
and hence we get the highest and lowest degrees bounded by $aw \pm
N/2$.  From this and the lemma the result follows.  $\Box$

\vspace{24pt}

{\it Proof (of Theorem 3).}
The skein definition of the Kauffman polynomial (Dubrovnik version) is given as
follows.
Let $D$ be the polynomial defined by the skein relation

\begin{eqnarray*}
 D_{K_+} - D_{K_-} &= & z ( D_{K_= } - D_{K_{)(} } ) \\
 D_{ K_{ < +>} } & = & \alpha D_{K} \\
 D_{ K_{< -> } } & = & \alpha^{-1} D_{K}
\end{eqnarray*}

where $K_+$, $K_-$ are link diagrams with positive and
negative crossings at a single particular crossing point
respectively
(and the rest of their diagrams coincide),
$K_= $ and $K _{)(}$ are link diagrams obtained by two ways
of smoothings at the crossing.
In the second and third equalities
$K_{< +>}$ (resp. $K_{< -> }$) denotes the link diagram
$K$ with a small positive (resp. negative) kink
added.

Now the Kauffman polynomial \cite{K4} is defined by

$$ Y_K(\alpha, z) = \alpha ^{-w(K)} D(\alpha, z) $$

where $w(K)$ is the writhe of the diagram $K$.

Note that the loop value is
$\mu =  z^{-1} (\alpha - \alpha^{-1}) +1 $.

As in the case of HOMFLY, the negative terms of $z$
come from the loop value, and do not contribute to
our estimates of the degrees in $t$.

Now we estimate the degrees with respect to $\alpha $, $\mu$, and $z$.
Each branch of the skein tree of $D_K$ (corresponding to
a state $\sigma$, which is a choice of smoothings of some of the
crossings that
results in a collection of framed unlinked unknots) will contribute
(a linear combination of) some  terms of $\alpha^p \mu^q z^r$.
Some of them may cancel out each other
but we only need an estimate of bounds of the degrees
so that we pick the highest and lowest possible degrees among them.

The degree of $z$ is the number of smoothings we performed
in the skein tree.
The degree of $\mu$ is the number of components of the state,
and the degree of $\alpha$ is the the writhe of the state.
Thus we obtain the following estimates:

$$
\begin{array}{rcccl}
d_{\max}(z) & \leq &\max\{n(\sigma)\}&\leq & c  \\
d_{\max}(\alpha) & \leq  & \max\{w(\sigma) \}  &
\leq & c_+\\
d_{\min}(\alpha) & \geq  & \min\{w(\sigma) \}  &
\geq & c_-\\
d_{\max}(\mu) & \leq & \max\{\ell(\sigma)\} -1&\leq &c +g-1 \\
\end{array}
$$
where $w(\sigma)$ and $\ell(\sigma)$  are the writhe and number of
components of $\sigma$ respectively,  $n(\sigma)$ is the number of
crossings smoothed to get to $\sigma$, and $c_+$, $c_-$, and $g$ are the
number of positive crossings, the number of negative crossings, and
the number of connected components of the diagram of $K$ respectively.

Thus after substitution $z=t^{1/2}-t^{-1/2}$ and  $\alpha = t^a$,
we get
\begin{eqnarray*}
d_{\max}(t) & \leq &  a \cdot d_{\max}(\alpha) + d_{\max}(z)/2 + (a
-1/2)d_{\max}(\mu) \\
&\leq& ac_+ + c/2 + (a-1/2)(c+g-1)\\
d_{\min}(t) & \geq &  a \cdot d_{\min}(\alpha) - d_{\max}(z)/2 - (a
-1/2)d_{\max}(\mu)\\
& \geq & ac_- - c/2 - (a-1/2)(c+g-1).
\end{eqnarray*}

The result then follows from the Lemma.
This gives a bound of
$$
N  \leq   3a c +(2a-1)(g-1),$$
but the same argument as in the proof of Theorem 1 gives
$$ N \leq 3ac. \Box$$

\begin{remark}
{\rm
It has been  conjectured that the Jones polynomial detects
knotting, {\it i.e.}, every nontrivial knot has a nontrivial
Jones polynomial.
By the above Corollary this can be restated as:
If a knot is nontrivial, then it has a nontrivial
Jones-Vassiliev invariant among those up to $N$,
the number of crossings of the given knot diagram.
     Moreover, we can consider the question whether there exist knots
whose knottedness is undetectable by {\em any} Vassiliev invariants of
small order in relation to the size of the knot diagram. Here ``small''
means an order less than or equal to the number of crossings in a
minimal diagram representing the knot, but it may be that one wants to
adjust the concept of small for this more general problem. It may be
useful to consider this generalization of the problem of knot detection
in the light of our work.
 }
\end{remark}

\end{document}